\documentclass[lettersize,journal]{IEEEtran}
\usepackage{amsmath,amsfonts,amsthm}
\usepackage{algorithmic}
\usepackage{algorithm}
\usepackage{array}
\usepackage{textcomp}
\usepackage{stfloats}
\usepackage{url}
\usepackage{verbatim}
\usepackage{graphicx,subfigure}
\usepackage{cite}
\usepackage{enumerate}
\usepackage{bm}
\usepackage{float}
\usepackage{booktabs} 
\usepackage{multirow}
\usepackage{mathrsfs}
\usepackage{comment}
\usepackage{xcolor}

\usepackage{fancyhdr}
\fancypagestyle{plain}{
  \fancyhf{} 
  \fancyhead[L]{\textcolor{blue}{This work has been submitted to the IEEE for possible publication. Copyright may be transferred without notice, after which this version may no longer be accessible.}} 
}

\theoremstyle{remark}
\newtheorem{myAssum}{Assumption}
\newtheorem{myRemark}{Remark}

\theoremstyle{plain}
\newtheorem{myTheo}{Theorem}

\begin{document}

\title{Ensemble Control for Stochastic Systems with Asymmetric Laplace Noises}

\author{Yajie Yu, Xuehui Ma, Shiliang Zhang, \IEEEmembership{Member, IEEE}, Zhuzhu Wang, Xubing Shi, \\Yushuai Li, \IEEEmembership{Member, IEEE}, Tingwen Huang, \IEEEmembership{Fellow, IEEE}
\thanks{This work was partially supported via the grant "Henan Province Science and Technology Tackling Key Problems" through grant agreement no. 232102210116. Yajie Yu is with Henan University of science and technology, Luoyang, China (yuyajie@haust.edu.cn). Xuehui Ma is with Xi'an University of Technology, Xi'an, China (xuehui.yx@gmail.com). Shiliang Zhang is with University of Oslo, Oslo, Norway (shilianz@ifi.uio.no). Zhuzhu Wang is with Stevens Institute of Technology, Hoboken, USA (zwang326@stevens.edu). Xubing Shi is with China Huaneng Group Co., Ltd., Zhengzhou,China. (sxb-hnhnbranch@163.com). Yushuai Li is with Aalborg University, Aalborg, Denmark (yushuaili@ieee.org). Tingwen Huang is with Texas A\&M University at Qatar, Doha 23874, Qatar (tingwen.huang@qatar.tamu.edu).}
}



\maketitle

\thispagestyle{plain} 

\begin{abstract}
This paper presents an adaptive ensemble control for stochastic systems subject to asymmetric noises and outliers. Asymmetric noises skew system observations, and outliers with large amplitude deteriorate the observations even further. Such disturbances induce poor system estimation and degraded stochastic system control. In this work, we model the asymmetric noises and outliers by mixed asymmetric Laplace distributions (ALDs), and propose an optimal control for stochastic systems with mixed ALD noises. Particularly, we segregate the system disturbed by mixed ALD noises into subsystems, each of which is subject to a specific ALD noise. For each subsystem, we design an iterative quantile filter (IQF) to estimate the system parameters using system observations. With the estimated parameters by IQF, we derive the certainty equivalence (CE) control law for each subsystem. Then we use the Bayesian approach to ensemble the subsystem CE controllers, with each of the controllers weighted by their posterior probability. We finalize our control law as the weighted sum of the control signals by the sub-system CE controllers. To demonstrate our approach, we conduct numerical simulations and Monte Carlo analyses. The results show improved tracking performance by our approach for skew noises and its robustness to outliers, compared with single ALD based and RLS-based control policy.\\

\textit{Note to Practitioners}—Most real-world systems are subject to stochasticity and disturbances from diverse sources, \textit{e.g.}, varying/disturbing environments, device aging, and inaccurate/corrupted measurements. Deriving a stochastic control law with disturbance resistance necessitates learning the unknown system dynamics and modeling the disturbance. This work proposes an adaptive ensemble control that can benefit a wide range of control applications, especially those with unknown parameters and suffering from noises and outliers of unknown distribution. We achieve our control design via the integration of modeling noises and outliers and learning unknown system parameters in our control law derivation. Theoretical analyses and simulations demonstrate the merits of the proposed control with robustness against noises and outliers and a resilient tracking control performance. We anticipate our results to contribute to control practices vulnerable to disturbances from the external environment and inner system, \textit{e.g.}, autonomous vehicle control, building automation, and renewable integration in power grids.
\end{abstract}


\begin{IEEEkeywords}
Adaptive control, asymmetric Laplace distribution, outliers, quantile regression, stochastic systems.
\end{IEEEkeywords}

\section{Introduction}
Adaptive control for stochastic systems has drawn significant attention from academia and industry with its advantages in real-world control applications. The deriving of an adaptive controller uses a feedback mode, and the control performance heavily depends on the estimation of unknown parameters in the stochastic system~\cite{aastrom2008adaptive,filatov2004adaptive}. Parameter estimation for the stochastic system is susceptible to the noises in the system observations ~\cite{MA2022157,alessandri2016moving,9925228,9502407}, like Gaussian noises, asymmetric noises and outliers induced by, \textit{e.g.}, device malfunctions, measurement noises, transmission error, or adversary attack~\cite{ma2022active,9189668} \textit{etc.} Therefore, the resilience against noises and outliers is critical to ensure reliable and precise system parameters estimation and a successful stochastic system control~\cite{DBLP:journals/tnn/ZhangCYZH18,DBLP:journals/cssc/ZhangCYZH19,DBLP:journals/tac/MaZLQSH24,DBLP:conf/eucc/MaCZLQS24}.

Solutions have been proposed to handle noises in adaptive control. Kalman filtering (KF) and recursive least square (RLS)~\cite{khodarahmi2023review} are well-known methods for optimal parameter estimation and have been successfully utilized in adaptive control systems. KF and RLS assume the noises in the system are of Gaussian distribution, nevertheless, such an assumption hardly holds in practical control systems~\cite{MA2022157,alessandri2016moving,9925228,9456096,maAdaptiveQuantileControl2022,9782410}. \textit{E.g.}, the noises in economic, biological, or flight control systems can be of asymmetric characters or outliers that do not follow Gaussian distribution~\cite{altunbacs2019impact,huang2017quantile,8453641,9843908,7112463}. As a result, the parameter estimation by KF or RLS can be inefficient for non-Gaussian noises like asymmetric noises and outliers, and the adaptive control based on KF or RLS might be reduced due to the inaccurate system parameter estimation~\cite{kral2009functional,MA2022157}. 


Outlier attenuation for stochastic system control has been investigated to either detect and remove outliers or to improve robustness for system parameter estimation against outliers during the control~\cite{ferdowsi_online_2014,kral2009functional,domanski2022outliers}. \textit{E.g.}, Alessandri \textit{et al.}~\cite{alessandri2016moving} developed a robust approach for the system estimation of discrete-time linear systems subject to outliers. Zaccarian \textit{et al.}~\cite{alessandri2018stubborn} introduced saturated output injection for the observer design of linear systems whose measurements are exposed to outliers. Ma \textit{et al.} constructed a recursive filter for nonlinear stochastic systems subject to outliers in system measurements~\cite{ma_probability-guaranteed_2021}. Wang \textit{et al.} proposed a full search robust estimation to numerically estimate the unknown system parameters and outliers simultaneously~\cite{wang_location_2021}. While those outlier attenuation approaches provide resistance for stochastic system control against outliers, we emphasize that there is no clear boundary between outliers and noises in practical systems. In that, outlier attenuation in system control should be combined with the analyzing and modeling of noises, and due attention should be paid to asymmetric noises that cannot be handled by established approaches like KF and RLS.

An alternative to model asymmetric noises and outliers in stochastic systems is to use the asymmetric Laplace distribution (ALD). Differing from Gaussian distribution, ALD probability density function can be configured with peak, thick-tail, and skewness, making it an efficient way to represent symmetric and asymmetric noises and outliers~\cite{kozubowski2000multivariate,YU2001437,yang2016posterior}. Several control strategies integrating symmetric Laplace noises in the system have been investigated. \textit{E.g.}, state estimation for a linear state-space model with symmetric Laplace noise has been studied in~\cite{9411735}, where an approximate inference algorithm is presented and successfully applied for the tracking control of audio frequency. Kalman-filter-liked adaptive filters and robust filters have been proposed to estimate system states for linear and nonlinear stochastic systems corrupted by symmetric Laplacian noise~\cite{khawsithiwong2011adaptive,9326329,liu2021nonlinear}. However, these mentioned works do not consider the skewness - or the asymmetric features - when constructing ALD noises. A few works consider asymmetric characters in system state estimation. \textit{E.g.}, Andersson \textit{et al.}~\cite{andersson2020optimum} proposed an online Bayesian filtering method for time series models corrupted by ALD noises; Xu \textit{et al.}~\cite{xu2021robust} developed an expectation maximization algorithm based system identification for parameter varying linear systems. While these algorithms demonstrate improved parameter estimation for systems subject to ALD noises, they have not been applied and contributed to stochastic system control. Ma \textit{et al.}~\cite{maAdaptiveQuantileControl2022} was the first to develop an adaptive control for stochastic systems corrupted by ALD noise, where they designed a Bayesian quantile sum estimator (BQSE) to identify unknown system parameters. While Ma \textit{et al.}'s work considers single ALD noise in the stochastic system, the noises in practice are more than often complex and may not be described only by a single probability distribution~\cite{6654117}.

In this work, we consider multiple ALD noises in stochastic systems, where we define the distribution of such noises as mixed asymmetric Laplace distribution (ALD), and propose an adaptive control method for stochastic systems subject to such disturbances. We contribute to this topic in that (\romannumeral1) we develop the mixed asymmetric Laplace distributions to model asymmetric noises and outliers in system observations (\romannumeral2) we partition the stochastic system into subsystems, and design an iterative quantile filter (IQF) to estimate parameters for each of the subsystems using the observations corrupted by asymmetric noises (\romannumeral3) we adopt the certainty equivalence (CE) principle to decouple the parameter estimation and system control. We formulate the optimal control problem for each of the subsystems, and derive the sub-optimal CE control law using the estimated parameter by IQF, and (\romannumeral4) we apply the Bayesian approach to ensemble individual subsystems, and derive the adaptive ensemble control law as the sum of sub-system CE control laws weighted by their Bayesian posterior probabilities. We demonstrate how our algorithm improves the control performance in numerical simulations and compare our approach with the ideal benchmark optimal control, the single ALD based control, and the RLS based adaptive control.

The remainder of this paper is as follows. Section~\ref{section2} formulates the optimal control problem for the stochastic systems with mixed asymmetric Laplace noises. Section~\ref{section3} elaborates the proposed adaptive ensemble control, where Section~\ref{subsection1} designs the iterative quantile filter for the parameter estimation of stochastic systems with ALD noises, Section~\ref{subsection2} derives the certainty equivalence control law, and Section~\ref{subsection3} derives the Bayesian based ensemble control law. We conduct Monte Carlo simulations and analyze the results in Section~\ref{section4}. Section~\ref{section5} concludes this paper.

\section{Problem statement} \label{section2}
This work considers a canonical linear, single-input, single-output stochastic system shown below:
\begin{equation}\label{system1}
	\begin{split}
		y(k + 1) = &{b_1}(k)u(k) +  \cdots  + {b_m}(k)u(k - m + 1) \\& + {a_1}(k)y(k)  +  \cdots + {a_n}(k)y(k - n + 1), 
	\end{split}
\end{equation}
\begin{equation}\label{system2}
    z(k) = y(k) + e(k), 
\end{equation}
where $k = 1, \cdots , N-1 $, $ {b_i}(k) $ and $ {a_i}(k) $ are unknown system parameters, $y(k)$ is the system output, $u(k)$ denotes the control input, $z(k)$ represents the measurement of the system, $e(k)$ is the measurement noise. We have the following assumptions for the described system:

\begin{myAssum}\label{assump1}
The system order $m$ and $n$ are known.
\end{myAssum}

\begin{myAssum}\label{assump2}
The initial state information $I_0=\{u(-1), \cdots, $ $ u(-m+1), y(0), \cdots, y(-n+1)\}$ is known.
\end{myAssum}

\begin{myAssum}\label{assump3}
The measurement noise $e(k)$ is considered as mixed asymmetric
Laplace noise, the distributions of which can be described by the sum of asymmetric Laplace distributions (ALDs) as shown below:
\begin{equation}\label{ALDS}
    p(e(k)) = \sum_{i=0}^{s} \gamma^{i} \mathcal{ALD}(e^{i}(k): \tau^{i}(k), \mu^{i}(k), \delta^{i}(k)),
\end{equation}
where the parameter $\gamma^{i}$ is unknown and it satisfies
\begin{equation}\label{}
    \sum_{i=0}^{s} \gamma^{i} = 1, \quad \gamma^{i}>0,
\end{equation}
and the variable $s$, skewness parameter $\tau^{i}(k)$, location parameter $\mu^{i}(k)$, and scale parameter $\delta^{i}(k)$ are known. The ALD probability density function is defined below \cite{yuBayesianQuantileRegression2001}:
\begin{equation}
	\begin{aligned}
		f_{ALD}^i(x)=\frac{\tau^i(1-\tau^i)}{\sigma^i}\left\{ \begin{array}{rc}
			e^{-(1-\tau^i)\frac{|x-\mu^i|}{\sigma^i}}, x<\mu^i\\
			e^{-\tau^i\frac{|x-\mu^i|}{\sigma^i}}, x\geq\mu^i
		\end{array}\right. 
	\end{aligned}
\end{equation}
where the skewness parameter $0 < \tau^i < 1 $ and the scale parameter $ \sigma^i > 0 $.
\end{myAssum}

The objective of this work is to design a control policy to force the system described in (\ref{system1}) and (\ref{system2}) to track the desired trajectory $y_r(k)$. We quantify this objective by the following cost function
\begin{equation} \label{costfunction}	
	J = E\left\{ \left. \sum_{k=0}^{N-1}{{\left[ {y(k + 1) - {y_r}(k + 1)} \right]}^2}  \right | {{I_k}}\right\},
\end{equation}
where $ E\left\{ {\cdot \left| {{I_k}} \right.} \right\} $ represents expectation given the information $ I_k $, and $I_k=\{u(k-1), \cdots, $ $ u(0), y(k), \cdots, y(1), I_0\}$. 
The control aims to derive a sequence of control signal $\{u(k)\}_{k=1}^{N-1}$ to minimize the cost function (\ref{costfunction}) subject to the system described in (\ref{system1}) and (\ref{system2}). We formulate the control problem as
\begin{equation} \label{problem}	
    \begin{aligned}
    (P) \quad &\min_{\{u(k)\}_{k=1}^{N-1}}   E\left\{ \left. \sum_{k=0}^{N-1}{{\left[ {y(k + 1) - {y_r}(k + 1)} \right]}^2}   \right | {{I_k}} \right\}, \\
    s.t. \quad & y(k + 1) = {b_1}(k)u(k) +  \cdots  + {b_m}(k)u(k - m + 1) \\
    & \quad \quad \quad   + {a_1}(k)y(k)  +  \cdots + {a_n}(k)y(k - n + 1), \\
    &  z(k) = y(k) + e(k),
    \end{aligned}
\end{equation}
where the measurement noise $e(k)$ is mixed asymmetric Laplace noise described by (\ref{ALDS}).

\begin{myRemark}
In problem $(P)$, the measurement noise is asymmetric Laplace noise rather than white Gaussian noise. In that, popular solutions like Recursive Least Square (RLS) or Kalman Filter (KF), which are derived based on Gaussian noise assumption, are ineffective in learning the system parameters during the control~\cite{kral2009functional,MA2022157}. In this work, we are motivated to design an iterative estimator with different quantiles, so as to deal with the ALD noise in the measurements, and improve the system tracking performance.
\end{myRemark}

\begin{myRemark}
Problem $(P)$ is an optimal control problem within a considered horizon. While there exists the uncertainty of the parameters of the stochastic system and the learning and identifying of those parameters is required, such parameter learning leads to the coupling between control and learning within the considered horizon~\cite{10055961}. The coupled learning and control renders nontrivial derivation of stochastic system control laws\cite{GUO1995435,MILOJEConvergence}. Such a situation necessitates (i) the decoupling of learning and control, and (ii) the reformulation of the control problem into a sub-optimal control problem with efficient problem-solving programming, so as to result in a practical control policy.
\end{myRemark}

\begin{myRemark}
As described in (\ref{ALDS}), the measurement noise $e(k)$ is the weighted sum of different ALD noises, where the weights are unknown. Compared with systems corrupted by a single ALD noise, the unknown aggregated ALD noise induces challenges toward the derivation of control law due to (i) the unknown aggregation settings and (ii) the unclear impact of the noise aggregation on the control law derivation. To address the challenges, this work aims at de-aggregating the mixed ALD noise, and developing the system control that can ensemble the control laws subject to individual de-aggregated ALD noise.
\end{myRemark}

\section{Controller design}\label{section3}

This section details the design of our ensemble control law for an uncertain system with mixed ALD noise. Section~\ref{subsection1} presents the designed iterative quantile filter. The iterative quantile filter de-aggregates the mixed ALD noise into weighted single ALD noises, and facilitates the online learning of parameters of the system subject to mixed ALD measurement noise. In Section~\ref{subsection2}, we apply the certainty equivalence principle to decouple the control and learning. This decoupling results in the solving of the optimal control problem with parameter uncertainties in the stochastic system, and enables the derivation of the certainty equivalence control law. We derive the ensemble control law in Section~\ref{subsection3}. Particularly, we partition the system into subsystems, each of which is subject to a single ALD noise, and then we derive the certainty equivalence control law for each of the subsystems. We finalize the adaptive ensemble control by aggregating the derived control laws for individual subsystems based on their Bayesian posterior probabilities.

\subsection{Iterative Quantile  Filter}\label{subsection1}

Classic system parameter learning approaches, \textit{e.g.}, least square or Kalman filter, assume the noise in the system follows Gaussian distribution. Based on this assumption, they derive recursive filtering for linear discrete-time systems by minimizing the posterior expected loss criterion. Nevertheless, least square or Kalman filter is based on mean regression that only works for white Gaussian noise. Mean regression merely considers the relationship between the mean of a response variable and its covariates, however, the mean of a response variable cannot describe the property of non-Gaussian noise, \textit{e.g.}, ALD noises~\cite{wangDistributedQuantileRegression2018,zhaoExpectationMaximizationApproach2016}. This section designs an iterative quantile filter (IQF) where we use quantile regression instead of mean regression approach. Compared with mean regression, quantile regression provides a more comprehensive description of the relationships between the response variable and its covariates by the estimation of a series of quantiles rather than only considering the mean. In that, quantile regression works for both Gaussian and non-Gaussian noises like ALD noises~\cite{wangDistributedQuantileRegression2018,zhaoExpectationMaximizationApproach2016}. In this section, we first brief how quantile regression works, and then detail the design of IQF in learning the parameters of the system subject to mix ALD noise.

Let $ Y $ be a random variable with cumulative distribution function $ {F_Y}(y) = P(Y \le y) $.
The $ \tau $-th $ (0<\tau<1)$ quantile of $ Y $ is given by $ {q_Y}(\tau ) = \inf\left\{ {y:{F_Y}(y) \ge \tau } \right\}$.
The parametric quantile regression can be expressed as below~\cite{wangDistributedQuantileRegression2018}:
\begin{equation}	
{{\hat y}_{(\tau )}} = \bm x^T \bm {\hat {\beta}}_{\tau},
\end{equation}
where $ {\hat y}_{(\tau )} $ represents the \(\tau\)-th quantile of the conditional distribution of the response variable \(y\), given the predictor variables \(\bm x\). $ {{\hat {\bm{\beta} } }_\tau } $ is the coefficient vector associated with the quantile \(\tau\).
$ {{\hat {\bm{\beta} } }_\tau } $ can be calculated by minimizing the non-differentiable sum of pinball-losses shown below:
\begin{equation} \label{loss1}
	\hat {\bm \beta} = \mathop {\min }\limits_{\bm \beta}  \sum {{\rho _\tau }} (y-\bm x^T \bm {\hat {\beta}}_{\tau}),
\end{equation}
where the loss function is $ {\rho _\tau }(u) = u\left( {\tau  - {\mathbb{I}_{(u < 0)}}} \right) $.
\(\mathbb{I}\) is the indicator function. Below we show how we design the iterative quantile filer based on quantile regression and illustrate the learning of unknown parameters of the system subject to mixed ALD noises using the designed filter.

We denote the unknown parameters in system (\ref{system1}) as a single vector,
\[\bm w(k) = [{b_1}(k), \cdots ,{b_m}(k),{a_1}(k), \cdots ,{a_n}(k)]^T.\]
The system (\ref{system1}) and (\ref{system2}) can be rewritten more concisely as 
\begin{equation}\label{system_CE}
	y(k + 1) = {\bm x}^T(k){\bm w}(k),
\end{equation}
and
\begin{equation}\label{system_IQF}
	z(k + 1) = {\bm x}^T(k){\bm w}(k)+ e(k + 1),
\end{equation}
where 
\begin{equation*}
	{\bm x}(k) = \left[ u(k), \ldots ,u(k - m+1), z(k), \ldots ,z(k -n+ 1) \right]^T.
\end{equation*}
The unknown parameters $\bm{w}$ are learned by quantile regression shown in \eqref{loss1}. We achieve the parameter learning by minimizing the non-differentiable sum of pinball-losses $\mathcal{L}_1({\bm{w}})$ below:
\begin{equation}
	\mathcal{L}_1({\bm{w}}) = E \sum\limits_{k = 1}^N {{\rho _\tau }\left[z(k+1)-{\bm x ^T}(k)\bm{w}(k)-e(k + 1)\right]}, 
\end{equation}
where  $ N $ is the number of observation samples. Relax the non-differentiable $L_1$-norm to the differentiable $L_2$-norm of the corresponding vector  \cite{maAdaptiveQuantileControl2022}, and the loss function is adjusted to the following:
\begin{equation}
	\begin{aligned}\label{loss2}
		&\mathcal{L}_2({\bm{w}}) = E \sum_{k=1}^{N}\left\{(1-\tau)\sum_{z(k+1)<{\bm x ^T}(k) \bm{w}(k)}\left[z(k+1) \right.\right.\\
		& \left. -{\bm x ^T}(k) \bm{w}(k)-e(k + 1)\right]^2+\tau\sum_{z(k+1)\geq{\bm x ^T}(k) \bm{w}(k)}\left[z(k+1)\right.\\
		&\left.\left.-{\bm x ^T}(k) \bm{w}(k)-e(k + 1)\right]^2 \right\},
	\end{aligned}
\end{equation}
which can be solved by recursive least squares. Theorem \ref{theory1} provides the full deviation of the IQF solution.

\begin{myTheo}\label{theory1}
Consider a stochastic linear system described by (\ref{system_IQF}), where the noise follows asymmetric Laplace distribution $\mathcal{ALD}(e(k): \tau(k), \mu(k), \delta(k))$. Initiate the system parameter as $\hat{w}(0)$ and estimation covariance as $\bm{P}(0)$. The parameter vector $\bm{w}(k)$ can be learned online by:
\begin{equation}
	\begin{aligned}
		&\bm{K}(k)=\frac{p(k)\bm{P}(k-1)\bm{x}(k)}{\bm{I}+p(k)\bm{x}^T(k)\bm{P}(k-1)\bm{x}(k)},\\
		&\hat{\bm{w}}(k+1)=\hat{\bm{w}}(k)+\bm{K}(k)\left[z(k+1)-{\bm x ^T}(k)\hat{\bm{w}}(k)-\varepsilon(k)\right],\\
		&\bm{P}(k+1)=[I-\bm{K}(k){\bm x ^T}(k)]\bm{P}(k),\\
        &\hat{\bm{w}}(0)=\bm{w}_0,\\
        &\bm{P}(0)=\bm{P}_0,
	\end{aligned}
\end{equation}
where the parameter $\varepsilon(k)$ is 
\begin{equation}\label{mean_ALD}
	\begin{aligned}
        \varepsilon(k) = \mu(k)+\frac{1-2\tau(k)}{\tau(k)(1-\tau(k))}\delta(k),\\
	\end{aligned}
\end{equation}
and the parameter $p(k)$ is 
\begin{equation}\label{para_p}
	\begin{aligned}
		p(k) =\left\{ \begin{array}{rc}
			1-\tau(k), \quad z(k+1)<{\bm x ^T}(k)\hat{\bm{w}}(k)\\
			\tau(k), \quad z(k+1)\geq{\bm x ^T}(k)\hat{\bm{w}}(k)
		\end{array}\right. .
	\end{aligned}
\end{equation}

\end{myTheo}

\begin{proof}
Define the weighted matrix $\bm{M}(k)$ as 
\begin{equation}\label{matrix_M}
	\begin{aligned}
		\bm{M}(k)= \left[ \begin{array}{cr}
			\bm{M}(k-1) & \mathbf{0}  \\
			\mathbf{0}  & p(k)
		\end{array} \right],
	\end{aligned}
\end{equation}
where the element $p(k)$ is set as (\ref{para_p}).  The vectors $\bm{X}(k)$, $\bm{Z}(k)$ and $\bm{E}(k)$ are defined as
\begin{equation}\label{vec_X}
	\begin{aligned}
		\bm{X}(k)= \left[ \begin{array}{c}
			\bm{X}(k-1) \\
			\bm{x}(k)
		\end{array} \right],
	\end{aligned}
\end{equation}
\begin{equation}\label{vec_Z}
	\begin{aligned}
		\bm{Z}(k)= \left[ \begin{array}{c}
			\bm{Z}(k-1) \\
			z(k)
		\end{array} \right],
	\end{aligned}
\end{equation}
\begin{equation}\label{vec_E}
	\begin{aligned}
		\bm{E}(k)= \left[ \begin{array}{c}
			\bm{E}(k-1) \\
			\varepsilon(k)
		\end{array} \right],
	\end{aligned}
\end{equation}
where the element $\varepsilon(k)$ is the mean of ALD noise $e(k)$ calculated by (\ref{mean_ALD}). With the definition (\ref{matrix_M}), (\ref{vec_X}), (\ref{vec_Z}) and (\ref{vec_E}), the loss function shown in (\ref{loss2}) could be rewritten as
\begin{equation}
	\begin{aligned}\label{loss3}
		\mathcal{L}_3[\bm{w}(k)] =& [\bm{Z}(k)-\bm{X}(k)\bm{w}(k)-\bm{E}(k)]^T\bm{M}(k)\\
                    &[\bm{Z}(k)-\bm{X}(k)\bm{w}(k)-\bm{E}(k)].
	\end{aligned}
\end{equation}
By setting the gradient of the loss (\ref{loss3}) to zero, we obtain the estimation of $\hat{\bm w}$ at the $k$-th instant as
\begin{equation}\label{w_hat}
	\begin{aligned}
		\hat{\bm w}(k) = [\bm{X}(k)^T\bm{M}(k)\bm{X}(k)]^{-1}\bm{X}(k)^T\bm{M}(k)[\bm{Z}(k)-\bm{E}(k)].
	\end{aligned}
\end{equation}
Define the matrix $\bm{P}(k)$ as
\begin{equation} \label{matrix_P}
	\begin{aligned}
		\bm{P}(k)&=[\bm{X}(k)^T\bm{M}(k)\bm{X}(k)]^{-1}.
	\end{aligned}
\end{equation}
Substituting (\ref{matrix_M}) and (\ref{vec_X}) into the matrix (\ref{matrix_P}), we obtain
\begin{equation} \label{matrix_P_inv}
	\begin{aligned}
		\bm{P}(k)&=[\bm{X}(k-1)^T\bm{M}(k-1)\bm{X}(k-1)+p(k)\bm{x}(k)\bm{x}^T(k)]^{-1}\\
		&=[\bm{P}^{-1}(k-1)+p(k)\bm{x}(k)\bm{x}^T(k)]^{-1}.
	\end{aligned}
\end{equation}
According to the matrix inversion lemma\cite{schott2016matrix}, matrix (\ref{matrix_P_inv}) can be rewritten as
\begin{equation}
	\begin{aligned} \label{p}
		\bm{P}(k)=\bm{P}(k-1)- \frac{p(k)\bm{P}(k-1)\bm{x}(k)\bm{x}^T(k)\bm{P}(k-1)}{\bm{I}+p(k)\bm{x}^T(k)\bm{P}(k-1)\bm{x}(k)}.
	\end{aligned}
\end{equation}
By substituting the matrix $\bm{P}(k)$ (\ref{matrix_P}) into parameter $\hat{\bm w}(k)$ (\ref{w_hat}), we get
\begin{equation}\label{w_hat_1}
		\hat{\bm{w}}(k)=\bm{P}(k)\bm{X}(k)^T\bm{M}(k)[\bm{Z}(k)-\bm{E}(k)].
\end{equation}
Substitute (\ref{matrix_M}), (\ref{vec_X}) and (\ref{vec_Z}) into (\ref{w_hat_1}), and the learning parameter $\hat{\bm w}(k)$ can be written as
\begin{equation}\label{w_hat_2}
    \begin{aligned}
		\hat{\bm{w}}(k)&=\bm{P}(k)\{\bm{X}^T(k-1)\bm{M}(k-1)[\bm{Z}(k-1)-\bm{E}(k-1)]\\
                        &+p(k)\bm{x}(k)[z(k)-\varepsilon(k)] \}.
    \end{aligned}
\end{equation}
From (\ref{w_hat_1}), we have 
\begin{equation}\label{}
    \begin{aligned}
		&\bm{X}^T(k-1)\bm{M}(k-1)[\bm{Z}(k-1)-\bm{E}(k-1)] \\
       & =\bm{P}^{-1}(k-1)\hat{\bm{w}}(k-1),
   \end{aligned}
\end{equation}
and substitute it into (\ref{w_hat_2}), we get
\begin{equation}\label{}
		\hat{\bm{w}}(k)=\bm{P}(k) \{ \bm{P}^{-1}(k-1)\hat{\bm{w}}(k-1)+p(k)\bm{x}(k)[z(k)-\varepsilon(k)] \}.
\end{equation}
We obtain the inverse of matrix $\bm{P}(k-1)$ by reforming the equation (\ref{matrix_P_inv})
\begin{equation} \label{matrix_P_inv_1}
		\bm{P}^{-1}(k-1)=\bm{P}^{-1}(k)-p(k)\bm{x}(k)\bm{x}^T(k).
\end{equation}
With (\ref{matrix_P_inv_1}), $\hat{\bm{w}}(k)$ could be written as
\begin{equation}\label{}
    \begin{aligned}
    \hat{\bm{w}}(k)=&\bm{P}(k)\{[\bm{P}^{-1}(k)-p(k)\bm{x}(k)\bm{x}^T(k)]\hat{\bm{w}}(k-1)\\
    &+p(k)\bm{x}(k)[z(k)-\varepsilon(k)]\}\\
    =&\hat{\bm{w}}(k-1)+p(k)\bm{P}(k)\bm{x}(k)\\
    &[z(k)-\bm{x}^T(k)\hat{\bm{w}}(k-1)-\varepsilon(k)].
    \end{aligned}
\end{equation}
We define the gain matrix $\bm{K}(k)$ as
\begin{equation} \label{K}
\bm{K}(k)=p(k)\bm{P}(k)\bm{x}(k),
\end{equation}
and then $\hat{\bm{w}}(k)$ is
\begin{equation}
	\hat{\bm{w}}(k) =\hat{\bm{w}}(k-1)+\bm{K}(k)[z(k)-\bm{x}^T(k)\hat{\bm{w}}(k-1)-\varepsilon(k)].
\end{equation}
Substituting equation (\ref{p}) into equation (\ref{K}) results in
\begin{equation} \label{K1}
	\bm{K}(k)=\frac{p(k)\bm{P}(k-1)\bm{x}(k)}{\bm{I}+p(k)\bm{x}^T(k)\bm{P}(k-1)\bm{x}(k)}.
\end{equation}
According to equation (\ref{p}) and (\ref{K1}), it can be obtained that
\begin{equation}
	\bm{P}(k)=[\bm{I}-\bm{K}(k)\bm{x}^T(k)]\bm{P}(k-1).
\end{equation}

\end{proof}

\begin{myRemark}
When $ \tau $ equals 0.5, the loss function $ \mathcal{L}_2({\bm{w}}) $ is coincided with the loss of the mean regression. Hence, the designed IQF can be regarded as a generalization of RLS or KF which only works for mean regression.
\end{myRemark}

\subsection{Certainty Equivalence Control}\label{subsection2}

Solving Problem (\ref{problem}) leads to a sequence of optimal control signals that can minimize the output tracking error in the whole control horizon. Nevertheless, the solving of Problem (\ref{problem}) requires that the system parameters are known, which does not hold in stochastic systems. The unknown system parameters also lead to coupling between system control and parameter learning, further complicating the control law deriving. To decouple the learning and control, we apply the certainty equivalence principle and facilitate a sub-optimal control law that can be efficiently derived and implemented.

According to the certainty equivalence principle, the parameters estimated by our IQF are viewed as true values when deriving the control law. Using the online estimated system parameters, we can obtain the cost function for a considered optimal control horizon. In this work, we reform the whole horizon optimal control costs to a one-step cost function to reduce the calculation load for control law. We define the cost function at the $k$-th instant as
\begin{equation} \label{loss3_new}	
	J_{CE}(k) = E\left\{ {{{\left[ {y(k + 1) - {y_r}(k + 1)} \right]}^2} \left| {I_k}  \right.} \right\}.
\end{equation}
Then we formulate the certainty equivalence control problem (P) for the uncertain system as
\begin{equation} \label{problem_sub}	
    \begin{aligned}
    (P_{CE}) \quad & \min_{u(k)|_{\bm{w}(k) = \hat{\bm{w}}(k)}} \quad J_{CE}(k),\\
    s.t. \quad \quad & y(k + 1) = {b_1}(k)u(k) +  \cdots  + {b_m}(k)u(k - m + 1) \\
    & \quad \quad \quad   + {a_1}(k)y(k)  +  \cdots + {a_n}(k)y(k - n + 1), \\
    &  z(k) = y(k) + e(k),
    \end{aligned}
\end{equation}
where the measurement noise $e(k)$ follows the distribution $\mathcal{ALD}(e(k): \tau(k), \mu(k), \delta(k))$. The solution to this certainty equivalence control problem is detailed in Theorem (\ref{theory2}).

\begin{myTheo}\label{theory2}
Consider the certainty equivalence control problem in $(P_{CE})$, the certainty equivalence control law is 
\begin{equation}
	\begin{aligned} \label{u_CE}
		u_{CE}(k)=\frac{y_r(k+1)-\bm{\eta}^T(k)\hat{\bm{\alpha}}(k)}{\hat{b}_1(k)},
	\end{aligned}
\end{equation}
where $ \hat{\bm{\alpha}}(k)=[\hat{b}_2(k), \cdots ,\hat{b}_m(k),\hat{a}_1(k), \cdots ,\hat{a}_n(k)]^T $, and $\bm{\eta}(k)=\left[ u(k-1), \ldots ,u(k - m+1), z(k), \ldots ,z(k -n+ 1) \right]^T$.

\end{myTheo}

\begin{proof}
Substitute equation (\ref{system_CE}) into equation (\ref{loss3_new}), and the cost function can be written as
\begin{equation}
J_{CE}(k) = E\left\{ {{{\left[ { {\bm x ^T}(k)\bm{w}(k) - {y_r}(k + 1)} \right]}^2} \left| {{I_k}} \right.} \right\}.
\end{equation}
With the parameter vector $ \bm{w}(k) $ partitioned as $ \bm{w}^T(k)=\left[ {b}_1(k) \quad \vdots \quad \bm{\alpha}^T(k)  \right] $ and the observation vector $ \bm{x}(k) $ partitioned as $ \bm{x}^T(k)=\left[ u(k) \quad \vdots \quad \bm{\eta}^T(k)  \right] $,
the cost function can be simplified as a function of $ u(k) $:
\begin{equation}\label{lu_preliminary}
	\begin{aligned}
		J_{CE}(k)=[b_1(k)u(k)+\bm{\eta}^T(k)\bm{\alpha}(k)-y_r(k+1)]^2.
	\end{aligned}
\end{equation}
The certainty equivalence control law can be derived assuming that the true values of unknown parameters are equal to the learned ones at each instant~\cite{boskovic2009certainty,karafyllis2018adaptive}. With this assumption, we have $ \bm{w}(k) =  \hat{\bm{w}}(k)=\left[\hat{b}_1(k) \quad \vdots \quad \hat{\bm{\alpha}}^T(k)\right] $, and the cost function can be rewritten as
\begin{equation}\label{lu}
	\begin{aligned}
		J_{CE}(k)=[\hat{b}_1(k)u(k)+\bm{\eta}^T(k)\hat{\bm{\alpha}}(k)-y_r(k+1)]^2.
	\end{aligned}
\end{equation}
The certainty equivalence control signal (\ref{u_CE}) can be obtained by setting the derivative of the equation (\ref{lu}) with respect to $ u(k) $ to zero.
\end{proof}
\begin{myRemark}
Since the measurement noise $e(k)$ is ALD noise, we use the IQF described in Theorem \ref{theory1} to learn the parameters $\hat{\bm{w}}$ iteratively. The CE control law is derive by $u_{CE}(k)=\left\{ \arg \min _{u(k)} J_{CE}(k) \right\} |_{\bm{w}(k) = \hat{\bm{w}}(k)}$.
\end{myRemark}

\begin{figure*}[tb]
\centering
\includegraphics[width=0.85\textwidth]{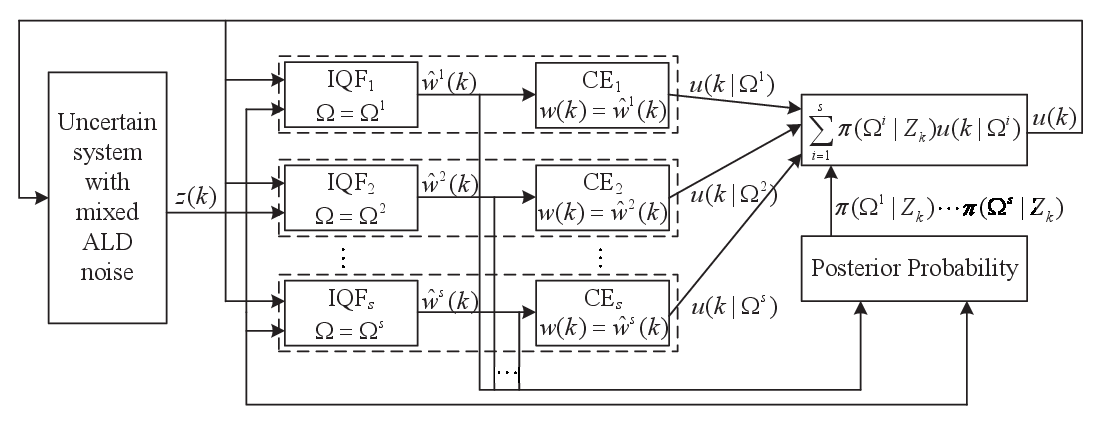}
\caption{The diagram of ensemble control law for an uncertain system with mixed ALD noises}
\label{fig1}
\end{figure*}

\subsection{Ensemble controller}\label{subsection3}

Fig. \ref{fig1} provides the block diagram of our ensemble control law for an uncertain system with mixed ALD noises. The ensemble control law $u^*(k)$ is achieved by aggregating control laws from subsystems each of which is associated with a single ALD noise. We use Bayesian method to gain the conditional posterior probability for each subsystem, and the conditional posterior probabilities serve as weights in control law aggregation. We detail the design of the ensemble controller below.

Theorem \ref{theory2} provides a certainty equivalence control law for uncertain systems corrupted by single ALD noises. This section extends this certainty equivalence control law with the capacity to handle mixed ALD noises.

According our assumption for mixed ALD noises in (\ref{ALDS}), we consider a system as an assembly of subsystems with different single ALD noises shown as $\mathcal{ALD}(e^{i}(k): \tau^{i}(k), \mu^{i}(k), \delta^{i}(k))$. Here we use a known finite set $\Omega$ to represent the subsystem with the parameters of $s$ ALD noises:
\begin{equation} \label{subsystem_def}	
    \begin{aligned}
        \Omega := \{(\tau^{1}, \mu^{1}, \delta^{1}), (\tau^{2}, \mu^{2}, \delta^{2}), \cdots, (\tau^{s}, \mu^{s}, \delta^{s})\}, 
    \end{aligned}
\end{equation}
and the $\Omega^i$ is the $i$-th subsystem with $i$-th ALD noise
\begin{equation} \label{subsystem_i}	
    \begin{aligned}
        \Omega^i = (\tau^{i}, \mu^{i}, \delta^{i}).
    \end{aligned}
\end{equation}
For each subsystem, we can learn its system parameters using IQF described in Theorem \ref{theory1}, and gain the certainty equivalence control laws by Theorem \ref{theory2}. Aggregating control laws for the subsystems results in the ensemble control law for mixed ALD noise. However, according to Assumption \ref{assump3}, the mixed weight $\gamma^i$ is unknown, thus hindering the control law aggregation. Here we apply the Bayesian method to calculate the posterior probabilities for each subsystem, and use the posterior probabilities as the aggregating weights for the ensemble control law. Theorem \ref{theory3} shows the derivation of the ensemble control law in detail. 

\begin{myTheo}\label{theory3}
Consider the control problem (\ref{problem}), and we aim to gain the ensemble control law below
\begin{equation}\label{u_ensemble}
	\begin{aligned}
		u^*(k)=\sum_{i=1}^{s}\pi(\Omega^i|Z_k)u(k|\Omega^i),
	\end{aligned}
\end{equation}
where $ \pi(\Omega^i|Z_k) $ is the discrete posterior probability computed by
\begin{equation}\label{bayesian_update}
	\begin{aligned}
		\pi(\Omega^i|Z_k)=\frac{f(z(k)|\Omega^i,Z_{k-1})\pi(\Omega^i|Z_{k-1})}{\sum_{j=1}^{s} f(z(k)|\Omega^j,Z_{k-1})\pi(\Omega^j|Z_{k-1})},
	\end{aligned}
\end{equation}
where $f(z(k)|\Omega^j,Z_{k-1})$ denotes the conditional probability density of $z(k)$ calculated by
\begin{equation}
    \begin{split}
	&f(z(k)|\Omega^i,Z_{k-1})=\frac{\tau^i(1-\tau^i)}{\sigma^i}\\
 &\exp\left\{ {-{\rho _{\tau^i} }(z(k)-{\bm x ^T}(k-1)\hat{\bm{w}}^i(k-1))} \right\},
    \end{split}
\end{equation}
and $Z_{k}$ is the measurement sequence $\{z(i)\}_{i=0}^k$.
$u(k|\Omega^i)$ is the certainty equivalence control law for the $i$-th subsystem with parameter $\Omega^i$, and is given by
\begin{equation}
	\begin{aligned} \label{u_CEi}
		u(k|\Omega^i)=\frac{y_r(k+1)-\bm{\eta}^T(k)\hat{\bm{\alpha}}^i(k)}{\hat{b}_1^i(k)},
	\end{aligned}
\end{equation}
where the parameter $\hat{\bm{w}}^i(k)=\left[\hat{b}_1^i(k) \quad \vdots \quad \hat{\bm{\alpha}}^{iT}(k)\right] $ is updated by 
\begin{equation}\label{param_learn_i}
	\begin{aligned}
      &\bm{K}^i(k)=\frac{p^i(k)\bm{P}^i(k-1)\bm{x}(k)}{\bm{I}+p^i(k)\bm{x}^T(k)\bm{P}^i(k-1)\bm{x}(k)},\\
		&\hat{\bm{w}}^i(k+1)=\hat{\bm{w}}^i(k)+\bm{K}^i(k)\left[z(k+1) \right.\\
 & \quad \quad \quad \quad \quad \left. -{\bm x ^T}(k)\hat{\bm{w}}^i(k)-\varepsilon^i(k) \right],\\
		&\bm{P}^i(k+1)=[I-\bm{K}^i(k){\bm x ^T}(k)]\bm{P}^i(k),\\
        &\hat{\bm{w}}^i(0)=\bm{w}_0,\\
        &\bm{P}^i(0)=\bm{P}_0,
	\end{aligned}
\end{equation}
where the parameter $\varepsilon^i(k)$ is 
\begin{equation}\label{mean_ALD_1}
	\begin{aligned}
        \varepsilon^i(k) = \mu^i(k)+\frac{1-2\tau^i(k)}{\tau^i(k)(1-\tau^i(k))}\delta^i(k),\\
	\end{aligned}
\end{equation}
and the parameter $p^i(k)$ is 
\begin{equation}
	\begin{aligned}
		p^i(k)=\left\{ \begin{array}{rc}
			1-\tau^i(k),\quad z(k+1)<{\bm x ^T}(k)\hat{\bm{w}}^i(k)\\
			\tau^i(k), \quad z(k+1)\geq{\bm x ^T}(k)\hat{\bm{w}}^i(k)
		\end{array}\right. .
	\end{aligned}
\end{equation}
\end{myTheo}

\begin{proof}
According to the certainty equivalence principle, the problem $(P)$ can be reformed into a certainty equivalence control problem described as
\begin{equation} \label{Problem_proof1}	
    \begin{aligned}
    (P_{CE}) & \min_{u(k)|_{\bm{w}(k) = \hat{\bm{w}}(k)}}  E\left\{\left[ y(k + 1) - {y_r}(k + 1) \right]^2 \left|I_k \right. \right\},\\
    s.t. \quad &y(k + 1) = {\bm x}^T(k){\bm w}(k), \\
         &z(k) = y(k) + e(k).
    \end{aligned}
\end{equation}
Solving problem $(P_{CE})$ is nontrivial since the measurement noise $e(k)$ is mixed ALD noises. With the definition of subsystem described in (\ref{subsystem_def}) and (\ref{subsystem_i}), we separate the problem $(P_{CE})$ into sub-problems each of which is associated with a single ALD noise $\Omega^i$. Using the smoothing property of conditional expectation, we can construct the cost function $J$ for the control problem in (\ref{Problem_proof1}) as
\begin{equation} \label{x}	
    \begin{aligned}
    \min J = \min_{u(k)} E\left\{ E\left\{ \left[ y(k + 1) - {y_r}(k + 1) \right]^2 \left|\Omega^i,I_k \right.\right\} \left|I_k \right. \right\}.
    \end{aligned}
\end{equation}
We interchange the minimization and the first expectation in (\ref{x}), resulting in
\begin{equation} \label{costf_appro}	
    \begin{aligned}
    \min J \approx E\left\{ \min_{u(k)} E\left\{ \left[ y(k + 1) - {y_r}(k + 1) \right]^2 \left|\Omega^i,I_k \right.\right\} \left|I_k \right. \right\}.
    \end{aligned}
\end{equation}
The inner minimization in (\ref{costf_appro}) represents the certainty equivalence control problem for the $i$-th subsystem described as
\begin{equation} \label{}	
    \begin{aligned}
    (P^i_{CE}) \quad & \min_{u(k)} \quad  E\left\{\left[ y(k + 1) - {y_r}(k + 1) \right]^2 \left|\Omega^i, I_k \right. \right\},\\
    s.t. \quad & y(k + 1) = {\bm x}^T(k){\bm w}(k), \\
    & z(k) = y(k) + e^i(k),
    \end{aligned}
\end{equation}
where the random noise $e^i(k)$ is the asymmetric Laplace noise $\mathcal{ALD}(e^{i}(k): \tau^{i}, \mu^{i}, \delta^{i})$. With the Theorem \ref{theory2}, we can obtain the certainty equivalence control law $u(k|\Omega^i)$ for the $i$-th subsystem, as is shown in (\ref{u_CEi}). The estimated system parameters $\hat{\bm{w}}^i(k)$ for calculating the control law $u(k|\Omega^i)$ can be obtained using Theorem \ref{theory1}. The learning of the $\hat{\bm{w}}^i(k)$ is shown as (\ref{param_learn_i}). To achieve the ensemble control, we use the posterior probability for subsystem with the $i$-th ALD noise $\Omega^i$ as the weight for aggregating control laws by subsystem. According to the Baye's rule, the posterior probability $\pi(\Omega^i|Z_k)$ follows
\begin{equation}
	\begin{aligned}
		\pi(\Omega^i|Z_k)&=\pi(\Omega^i|z(k),Z_{k-1}),
	\end{aligned}
\end{equation}
and the update of $\pi(\Omega^i|Z_k)$ can be achieved via equation (\ref{bayesian_update}). According to (\ref{costf_appro}), the outer expectation can be calculated as a probability weighted sum of inner values, where the probability weights are obtained as posterior probabilities $\pi(\Omega^i|Z_k)$ in (\ref{bayesian_update}). 
Our ensemble control law is calculated as a posterior probability weighted sum of subsystem control laws, where the $i$-th subsystem control law is given by the certainty equivalence control law (\ref{u_CEi}). The ensemble control law is finalized in the format shown in (\ref{u_ensemble}).
\end{proof}

\begin{myRemark}
The initial probability $\pi(\Omega^j|Z_0) $ is set as $1/s$, which means the initial probabilities of each subsystem are assumed to be equal.
\end{myRemark}

\section{Simulations and results} \label{section4}

\begin{figure*}[tb]\footnotesize	
	\begin{minipage}{0.333\linewidth}
		\vspace{0pt}        \centerline{\includegraphics[width=\textwidth]{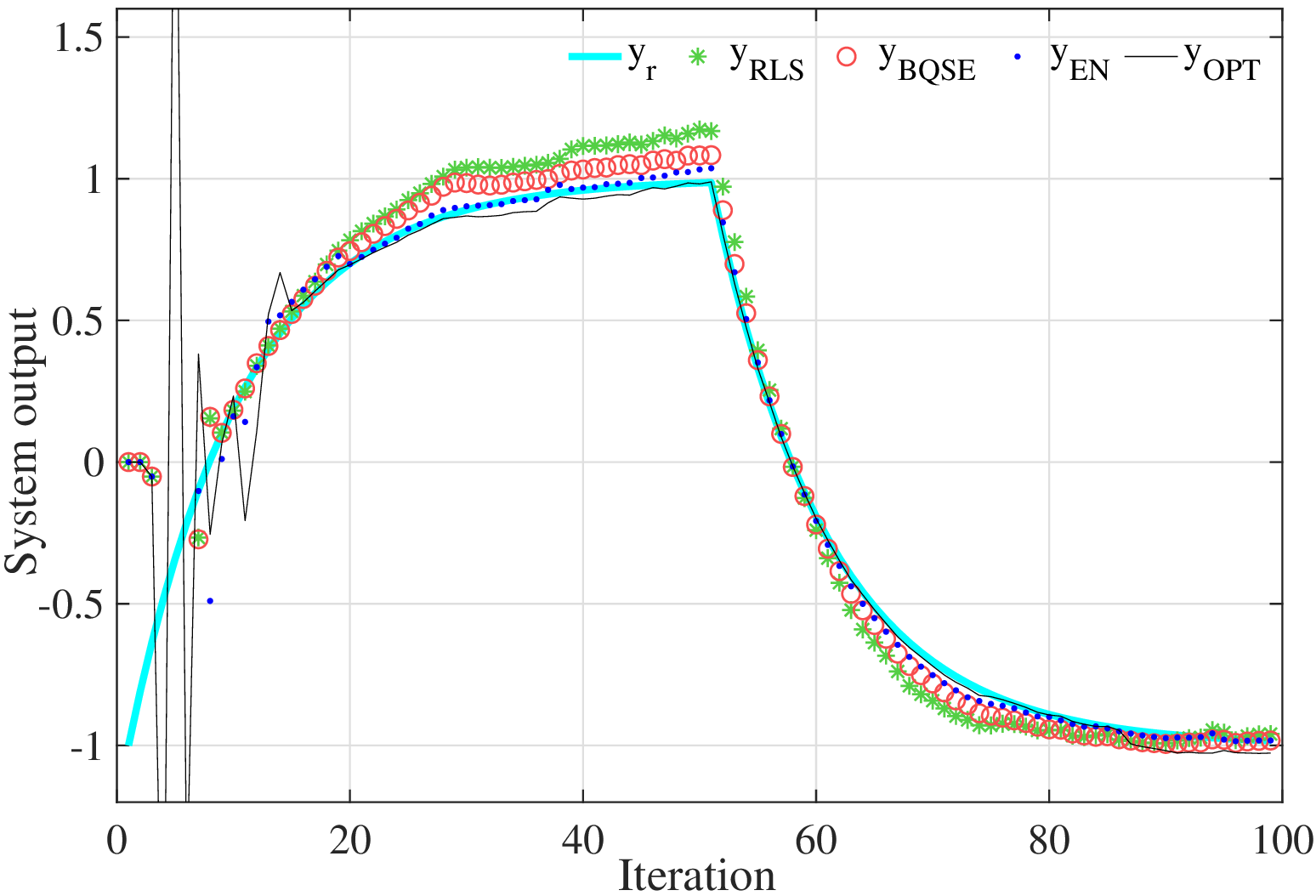}}
\centerline{(a) Square wave }
	\end{minipage}
	\begin{minipage}{0.333\linewidth}
		\vspace{0pt}		\centerline{\includegraphics[width=\textwidth]{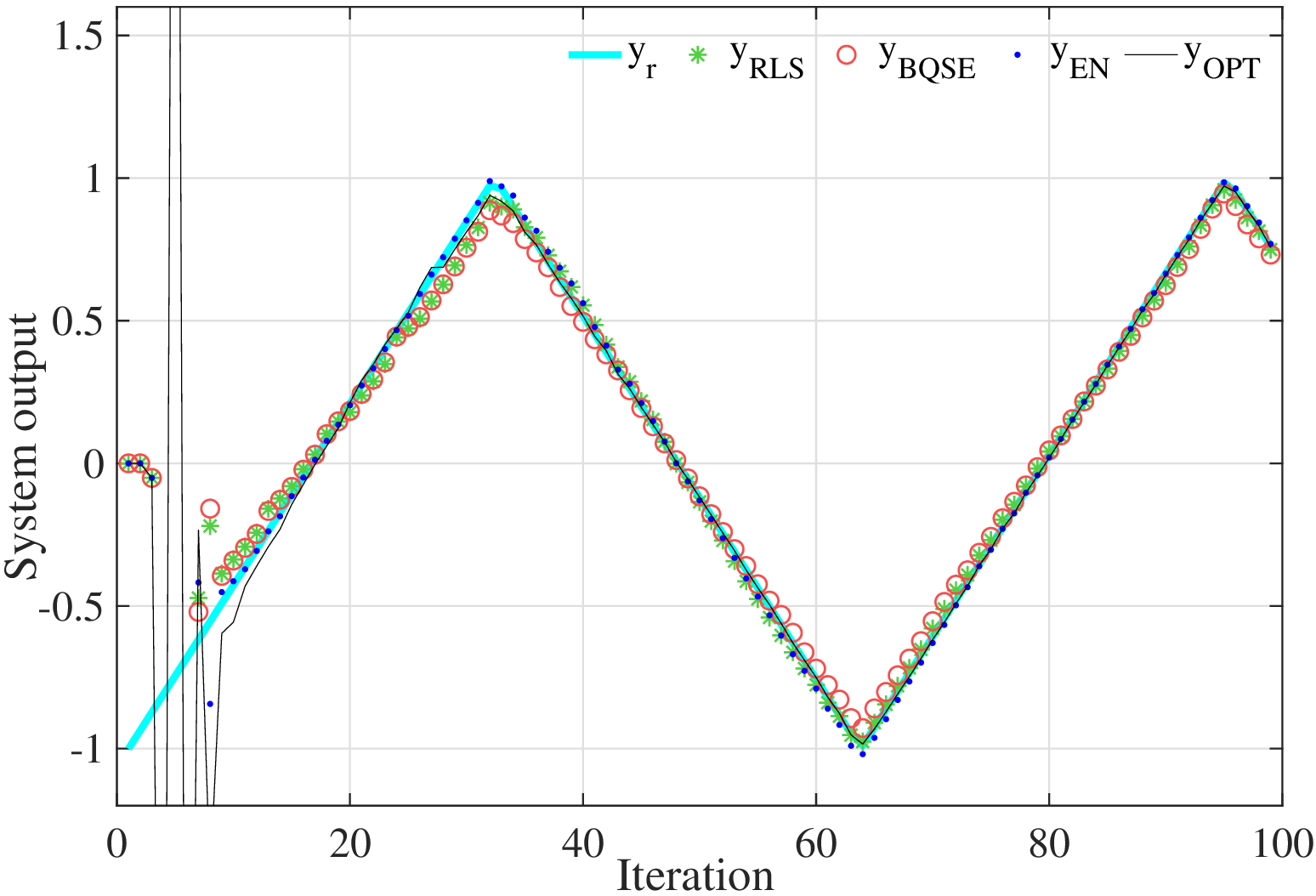}}	 
  \centerline{(b) Triangular wave}
	\end{minipage}
	\begin{minipage}{0.333\linewidth}
		\vspace{0pt}		\centerline{\includegraphics[width=\textwidth]{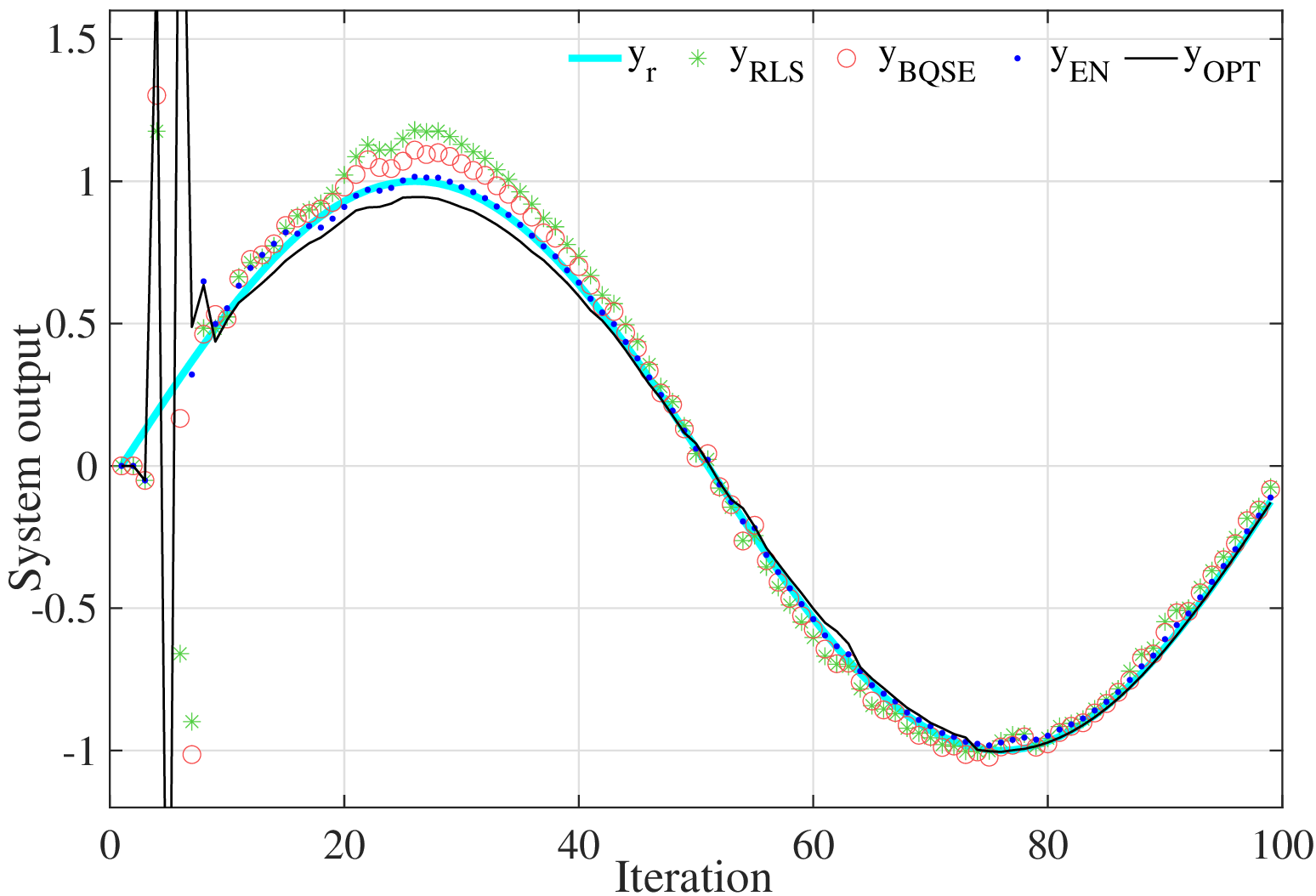}}	 
  \centerline{(c) Sine wave}
	\end{minipage} 
	\caption{The output of the system controlled by optimal control with known parameters, certainty equivalence control with RLS, and adaptive ensemble control, represented by $y_{OPT}$, $y_{RLS}$, and $y_{EN}$, respectively. The reference trajectory $y_r$ in (a), (b), and (c) are 0.01Hz square wave with unit amplitude filtered by the transfer function $1/(s+1)$, 0.01Hz triangular wave with 1 unit amplitude, and 0.01Hz sine wave with 1 unit amplitude, respectively.}
	\label{fig2}
\end{figure*}

This section conducts simulations implementing ensemble control on a minimum system. Algorithm \ref{algo1} summarizes the implementation of the designed ensemble control. We assess the tracking performance of the controller subject to mixed ALD noises under different tracking scenarios in Section \ref{exp_subs1}. Section \ref{exp_subs2} constructs different types of outliers and examines the tracking performance under various noise conditions and outliers.

\begin{algorithm}[htb] 
	\caption{Ensemble control law for an uncertain system with mixed ALD noises.} 
	\label{alg:Framwork} 
	\begin{algorithmic}[1] 
		\REQUIRE  
        The desired trajectory $y_r(k)$;    
		\STATE Initialize parameter vector $\hat{\bm{w}}^i(0)$, covariance $\bm{P}^i(0)$ and $\pi(\Omega^j|Z_0) $ for every subsystem, and state information $I_0$; 
        \STATE Set the parameters of each ALD noise $\Omega^i$;
        \STATE $k \leftarrow 1$ 
        \REPEAT
		\STATE Calculate the parameter vector for each subsystem $ \hat{\bm{w}}^i(k+1) $ according to equation (\ref{param_learn_i}); 
		\STATE Calculate the certainty equivalence control law for each subsystem $ u(k|\Omega^i) $ based on equation (\ref{u_CEi}); 
		\STATE Calculate the posterior probability for each subsystem $ \pi(\Omega^i|Z_k) $ according to equation (\ref{bayesian_update});
		\STATE Calculate the ensemble control law $ u^*(k) $ based on equation (\ref{u_ensemble});
		\STATE Apply $u^*(k)$ to the system; 
        \STATE $k \leftarrow k+1$ 
        \UNTIL $k \leftarrow N-1$ 
	\end{algorithmic} \label{algo1}
\end{algorithm} 

\subsection{Simulations of the control under mixed ALD noises } \label{exp_subs1}
This section considers the control of a minimum-phase system shown in equation (\ref{system1}) with system and ALD noise parameters below 
\begin{equation}
	\begin{aligned}
		a_1&=-1.41, \quad a_2=0.9, \quad n=2,\\
		b_0&=0.5, \quad m=0,\\
  p(e(k))& =  0.8 \mathcal{ALD}(e^{1}(k): 0.95, 0, 0.01)\\&+0.2 \mathcal{ALD}(e^{2}(k): 0.85, 0, 0.01).
	\end{aligned}
\end{equation}
We implement ensemble control on the system above and compare its tracking performance with the control by RLS, BQSE, and the ideal optimal controller where the coefficient $\gamma$ is known. This simulation uses three different reference trajectories for tracking: a square wave filtered by a network of transfer function $1/(s+1)  $, a triangle wave, and a sine wave with a unit amplitude at 0.1Hz.
The initial parameter $ {\bm w}^T(1) $ is set to $ [0.1,0.1,0.1] $, and $ {\bm P} = 100{\bm I} $. The prior probability $ \gamma(1) $  set as $1/s$, that is $ [0.5,0.5] $. We use same reference trajectories, noise sequence and initial conditions for the ensemble control, RLS, and BQSE.

Fig. \ref{fig2} shows the system output by RLS, BQSE, the ideal optimal control (OPT), and the designed ensemble control (EN) under various inference input conditions disturbed by mixed ALD noises. For all the three reference trajectories, we observe significant oscillations at the beginning of control across all the considered methods. Such oscillations disappear after the first 10 steps, indicating that the parameter learning works during the start-up of the control. We also observe from Fig.~\ref{fig2} that the outputs by OPT and EN track the reference signals closely, while RLS and BQSE well track the triangular wave trajectory yet diverge from the sine and square trajectories to different levels.

We conduct Monte Carlo simulations to evaluate the control performances and examine the statistics. We define $\bar J(i)$ as the accumulated error for the $i$-th Monte Carlo simulation. We define the control performance index as the average accumulated cost $\bar J$ below
\begin{equation}
\bar J = \frac{1}{{{N_{MC}}}}\sum_{i=1}^{N_{MC}} \bar J(i),
\end{equation}
where $ N_{MC} $ is the number of Monte Carlo simulations. We set $ N_{MC} $ as 100 in our simulation.

\begin{figure*}[htb]\footnotesize	
	\begin{minipage}{0.333\linewidth}
		\vspace{0pt}        \centerline{\includegraphics[width=\textwidth]{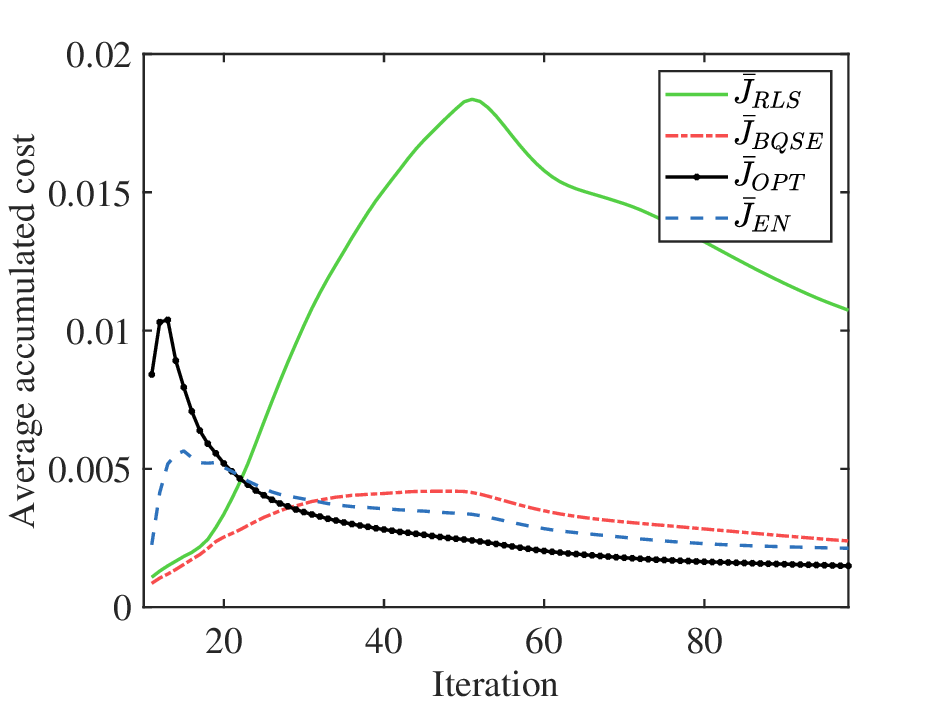}}
  \centerline{(a) Square  wave}
	\end{minipage}
	\begin{minipage}{0.333\linewidth}
		\vspace{0pt}		\centerline{\includegraphics[width=\textwidth]{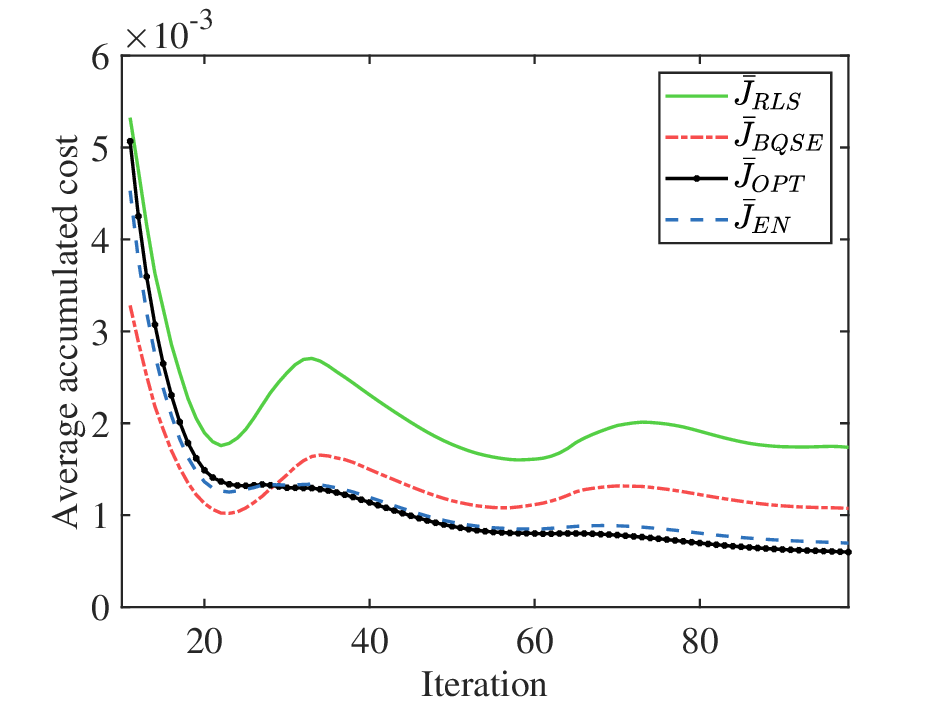}}	 
  \centerline{(b) Triangular wave}
	\end{minipage}
	\begin{minipage}{0.333\linewidth}
		\vspace{0pt}		\centerline{\includegraphics[width=\textwidth]{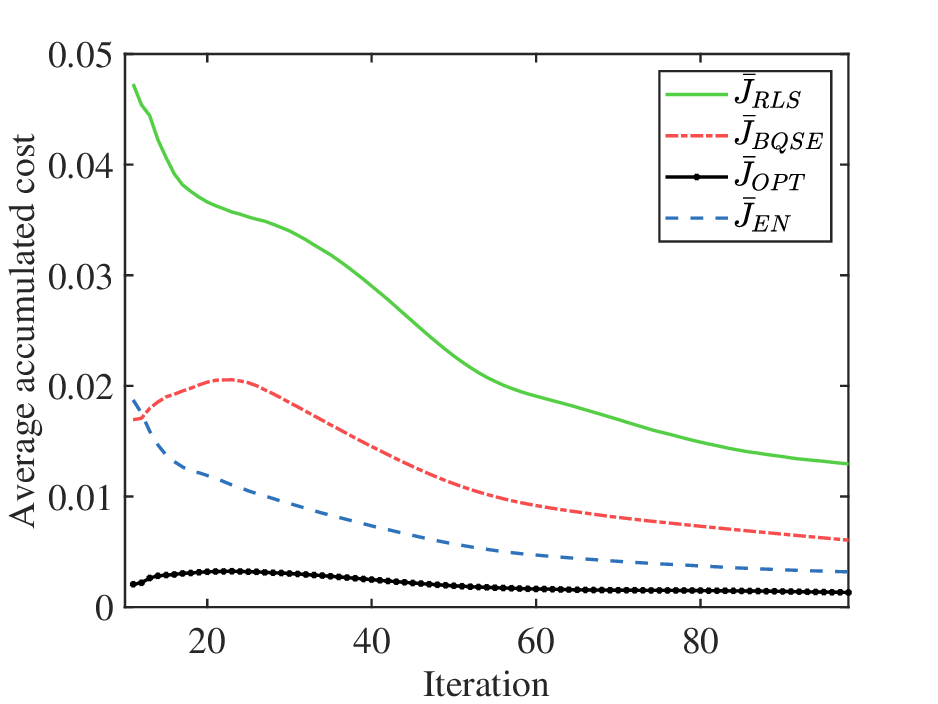}}	 
  \centerline{(c) Sine wave}
	\end{minipage} 
	\caption{The average accumulated error for 100 Monte Carlo simulations}
	\label{fig4}
\end{figure*}

\begin{figure*}[htb]\footnotesize	
	\begin{minipage}{0.333\linewidth}
		\vspace{0pt}        \centerline{\includegraphics[width=\textwidth]{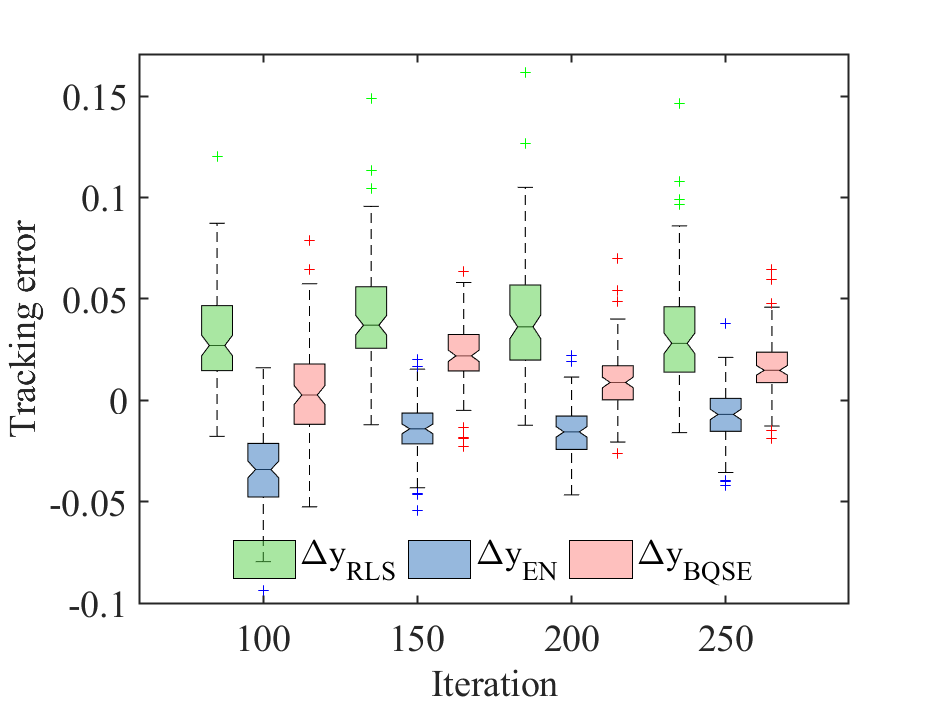}}
\centerline{(a) Square wave }
	\end{minipage}
	\begin{minipage}{0.333\linewidth}
		\vspace{0pt}		\centerline {\includegraphics[width=\textwidth]{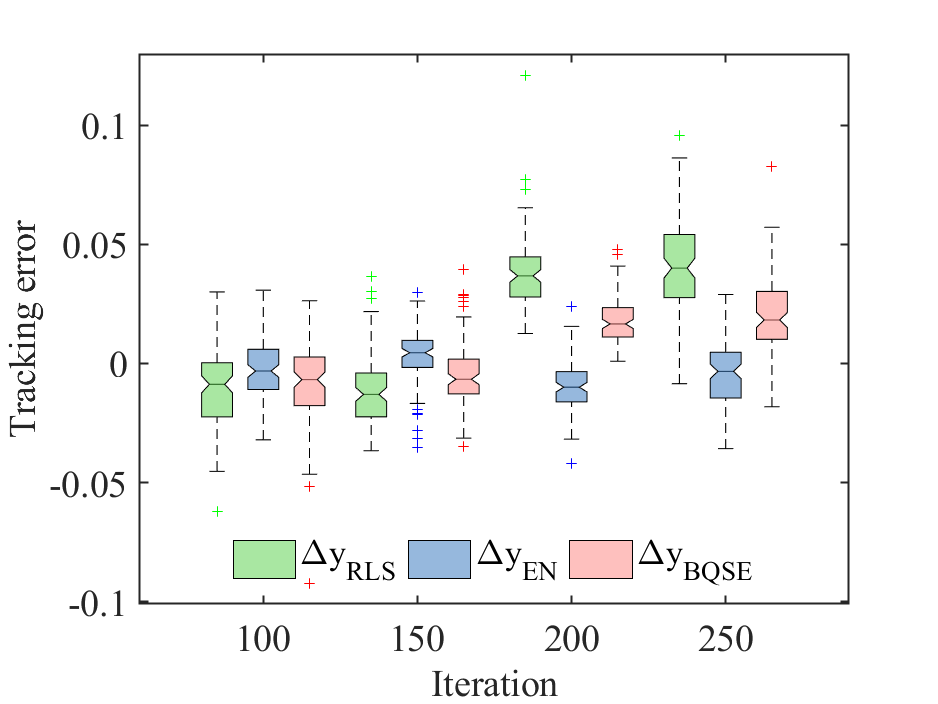}}
  \centerline{(b) Triangular wave}
	\end{minipage}
	\begin{minipage}{0.333\linewidth}
		\vspace{0pt}		\centerline{\includegraphics[width=\textwidth]{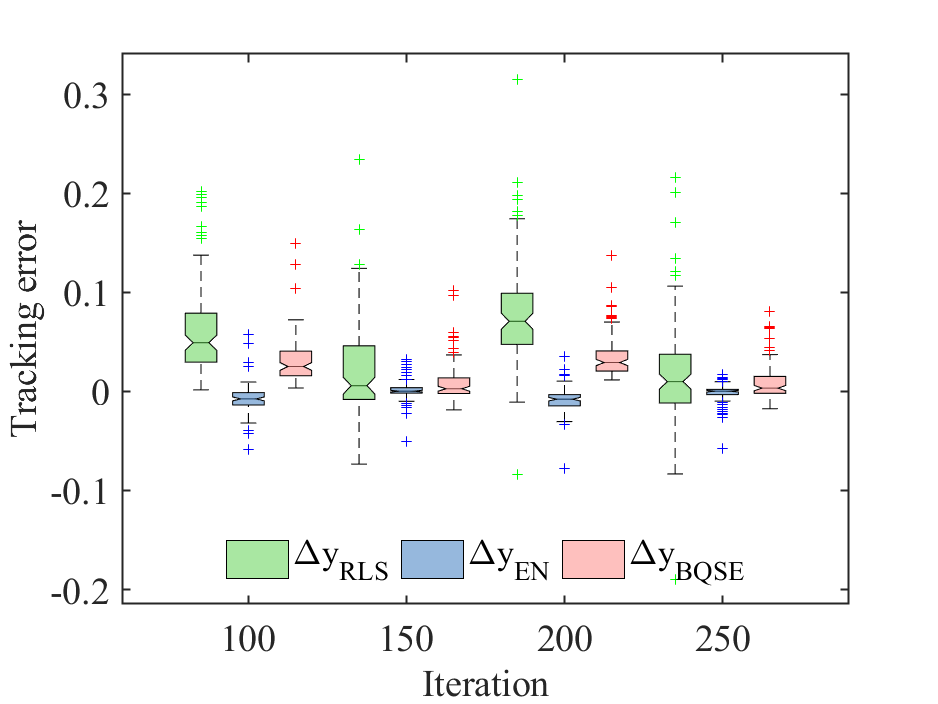}}	 
  \centerline{(c) Sine wave}
	\end{minipage} 
	\caption{Tracking error for the control of 100-300th iterations based
on 100 Monte Carlo simulations.}
	\label{fig5}
\end{figure*}

Fig. \ref{fig4} shows the tracking error from the 10th to 100th iteration for RLS, BQSE, OPT and EN based on 100 Monte Carlo simulations. For all the square, triangular, and sine trajectories, the average accumulated error curve for EN is the closest to that of OPT compared with BQSE and RLS, indicating the advantage of the designed ensemble control over BQSE and RLS. We summarize the average accumulated errors over the whole simulation horizon in Table \ref{tab1}. From Table \ref{tab1}, we observe that EN achieves the lowest average accumulated error compared with RLS and BQSE. Particularly, for sine trajectory, the average accumulated error by EN is 82.88\% lower than that of RLS and 20.83\% lower than that of BQSE. For triangular trajectory, the proposed ensemble control reduces the average accumulated error by 58.82\% and 36.36\% compared to RLS and BQSE, respectively. For square trajectory, the reduction of the average accumulated error by ensemble control is 74.80\%. and 46.55\% compared to RLS and BQSE, respectively.



\begin{table}[h]
	\caption{Average accumulated error by ensemble control, RLS, and BQSE from the 10th to 100th iteration.}
	\label{tab1}
	\centering      
	\begin{tabular}{p{2cm}p{1.5cm}p{1.5cm}p{1.5cm}}
		\hline
		\multirow{3}{*}{Methods} & \multicolumn{3}{c}{Reference trajectory}  \\
		\cmidrule(r){2-4}&Square&Triangle&Sin\\
		\hline
		RLS&0.0111&0.0017&0.0123\\
            \hline
		BQSE&0.0024&0.0011&0.0058\\
		\hline
		OPT&0.0016&0.0006&0.0013\\
		\hline
		Ensemble control&0.0019&0.0007&0.0031\\
		\hline
	\end{tabular}
\end{table}

\begin{table}[h]
	\caption{Average accumulated cost for different input references  from the 100th to 300th iteration}
	\label{tab2}
	\centering      
	\begin{tabular}{p{2cm}p{1.5cm}p{1.5cm}p{1.5cm}}
		\hline
		\multirow{3}{*}{Methods} & \multicolumn{3}{c}{Reference trajectory}  \\
		\cmidrule(r){2-4}&square&triangle&sin\\
		\hline
		RLS&0.0014&0.0011&0.0062\\
  		\hline		
		BQSE&0.0004&0.0003&0.0010\\
		\hline		
		Ensemble control&0.0003&0.0001&0.0003\\
		\hline
	\end{tabular}
\end{table}

We conduct further simulations to analyze the tracking performance of our ensemble control from 100th to 300th iteration, during which the controller has learned the system parameters well and thus is stable. Fig. \ref{fig5} shows the box plot of the tracking error of the ensemble control in comparison with RLS and BQSE, based on the results of 100 Monte Carlo simulations. We observe that the designed ensemble control provides more concentrated boxes than RLS and BQSE, especially for sine trajectory, indicating better and stable control performance by the proposed ensemble controller. Table \ref{tab2} summarizes the average accumulated error from the 100th to 300th iteration based on 100 Monte Carlo simulations. The average accumulated error by the proposed controller is significantly lower than that of RLS for all the three types of trajectories for the tracking, while the average accumulated error by the ensemble control is 25\% lower than that of BQSE for sine inputs, 66.67\% for triangular wave inputs, and 70.00\% for square wave inputs, respectively.
The simulation results demonstrate that the proposed ensemble controller is effective in tracking control subject to mixed ALD noises.

\subsection{Simulations with various noise conditions and outliers} \label{exp_subs2}

\begin{figure}[tb]	
	\centering	
        \includegraphics[width=0.48\textwidth]{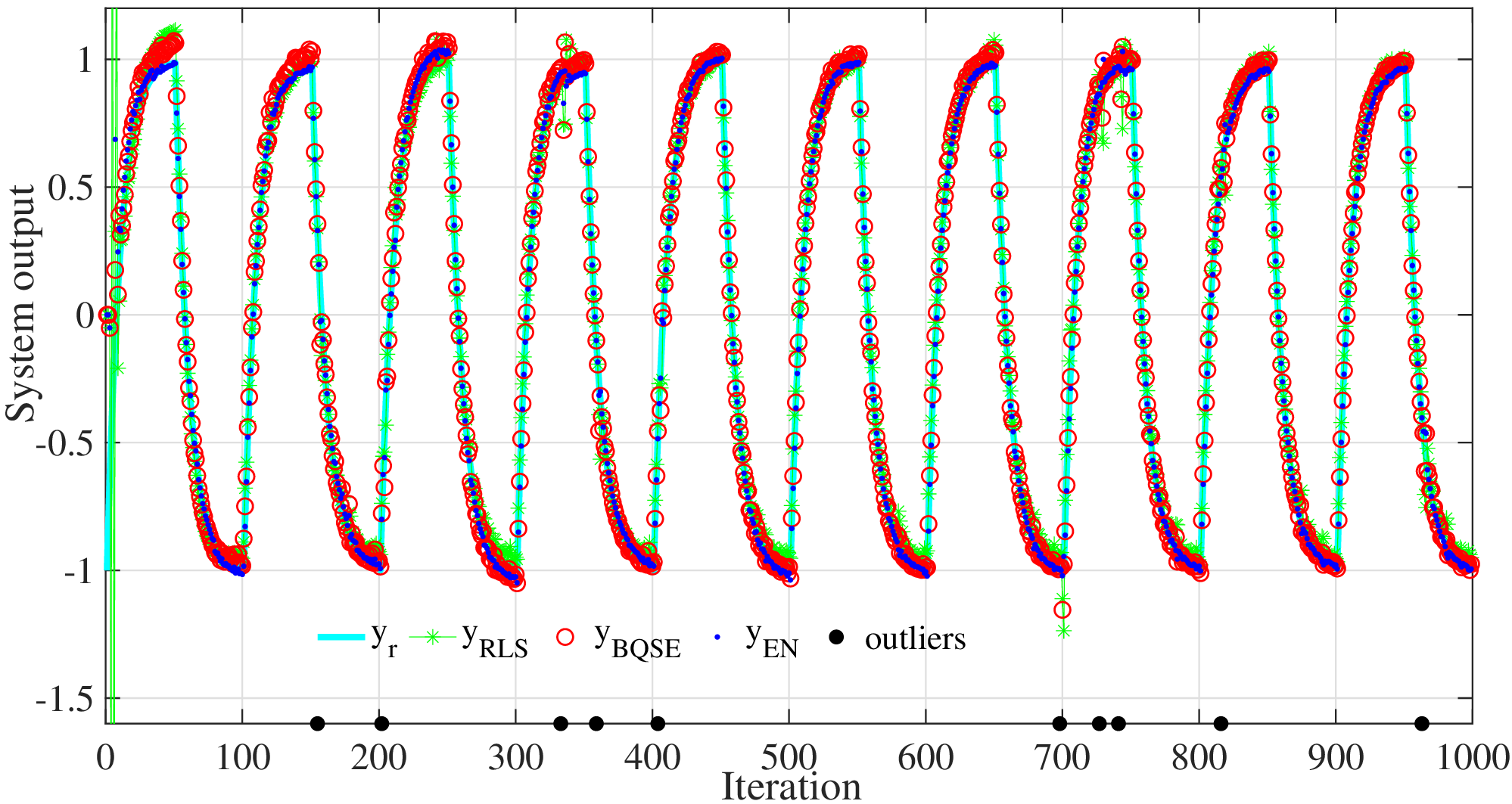}	
        \caption{The system output subject to Noise-\uppercase\expandafter{\romannumeral 1}}
	\label{fig6}	
\end{figure}

\begin{figure}[tb]	
	\centering	
        \includegraphics[width=0.48\textwidth]{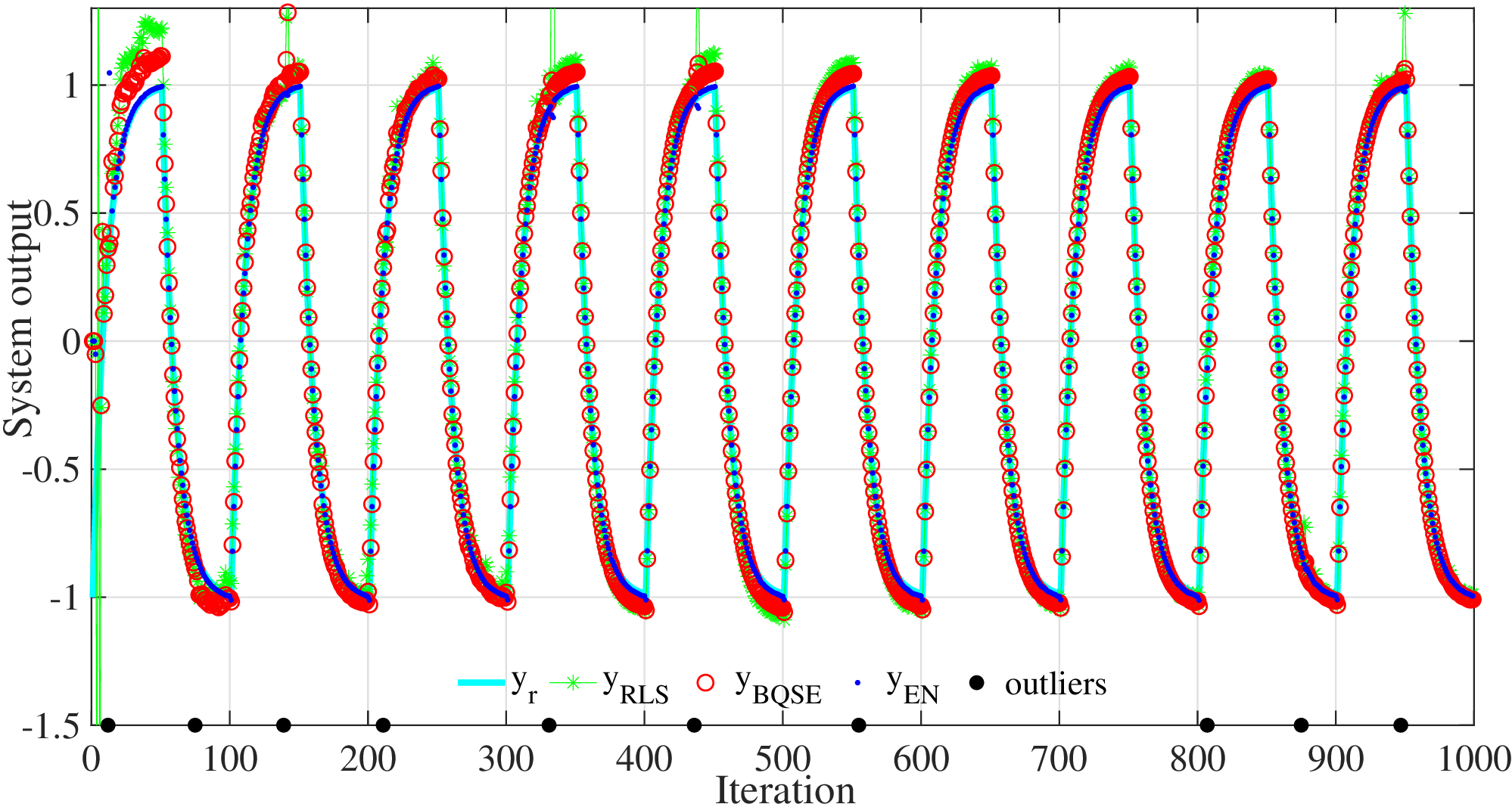}	
        \caption{The system output subject to Noise-\uppercase\expandafter{\romannumeral 2}}
	\label{fig7}	
\end{figure}

We construct various noise conditions and different types of outliers to test the robustness of the proposed ensemble control, especially its resilience against outliers. We use the following specifications when constructing mixed ALD noises and outliers:
\begin{enumerate}[{Noise}-I]
  \item \begin{equation*}
\begin{aligned}
    p(e(k)) =  &0.99 \mathcal{ALD}(e^{1}(k): 0.95, 0, 0.01)\\&+0.01 \mathcal{ALD}(e^{2}(k): 0.85, 2, 0.01),
\end{aligned}
\end{equation*}
  \item \begin{equation*}
\begin{aligned}
    p(e(k)) =  &0.99 \mathcal{ALD}(e^{1}(k): 0.95, 0, 0.01)\\&+0.01 \mathcal{ALD}(e^{2}(k): 0.85, 0, 2),
\end{aligned}
\end{equation*}
  \item \begin{equation*}
\begin{aligned}
    p(e(k)) =  &0.99 \mathcal{ALD}(e^{1}(k): 0.95, 0, 0.01)\\&+0.01 \mathcal{N}(e^{2}(k): 2, 0.01),
\end{aligned}
\end{equation*}
\item \begin{equation*}
\begin{aligned}
    p(e(k)) =  &0.99 \mathcal{ALD}(e^{1}(k): 0.95, 0, 0.01)\\&+0.01 \mathcal{N}(e^{2}(k): 0, 2).
\end{aligned}
\end{equation*}
\end{enumerate}
The first component of all the noises is a basic ALD noise, while the second component models outliers. Noise-\uppercase\expandafter{\romannumeral 1} and Noise-\uppercase\expandafter{\romannumeral 2} both contains ALD  outliers with probability 0.01, with $\tau = 0.85$,  $\mu = 1$ and $\sigma^2 = 0.01$, and with $\tau = 0.85$, $\mu = 0$ and $\sigma^2 = 2$, respectively. Noise-\uppercase\expandafter{\romannumeral 3} and Noise-\uppercase\expandafter{\romannumeral 4} both contains Gaussian outliers with probability 0.01, with  $\mu = 1$ and $\sigma^2 = 0.01$, and with $\mu = 0$ and $\sigma^2 = 2$, respectively. We use the same system described in Section~\ref{exp_subs1} to test the tracking performance in this simulation, and we use the square trajectory described in Section~\ref{exp_subs1} for tracking control.

Fig. \ref{fig6} shows the output of the system subject to Noise \uppercase\expandafter{\romannumeral 1} by RLS, BQSE and the ensemble control. We observe that due to the outliers at the 333rd, 698th, 727th and 941st iterations, the outputs by RLS and BQSE exhibit noticeable spikes. Fig. \ref{fig7} shows the system outputs by RLS, BQSE, and the ensemble control under Noise \uppercase\expandafter{\romannumeral 2}. Since Noise \uppercase\expandafter{\romannumeral 2} has a larger variance than Noise \uppercase\expandafter{\romannumeral 1}, the outliers generated by Noise \uppercase\expandafter{\romannumeral 2} tends to have higher magnitudes. As a result, we observe significant spikes during the control in the presence of outliers at the 331st iterations for RLS and BQSE. Additionally, we observe at the 139th, 436th and 947th iterations that, while the RLS output exhibits noticeable spikes, the BQSE output shows relatively small spikes. In contrast to RLS and BQSE, the system output by the proposed ensemble control effectively tracks the reference trajectory under both Noise-\uppercase\expandafter{\romannumeral 1} and Noise-\uppercase\expandafter{\romannumeral 2}, manifesting its resilience of control against outliers.

\begin{figure}[h]	
	\centering	
        \includegraphics[width=0.48\textwidth]{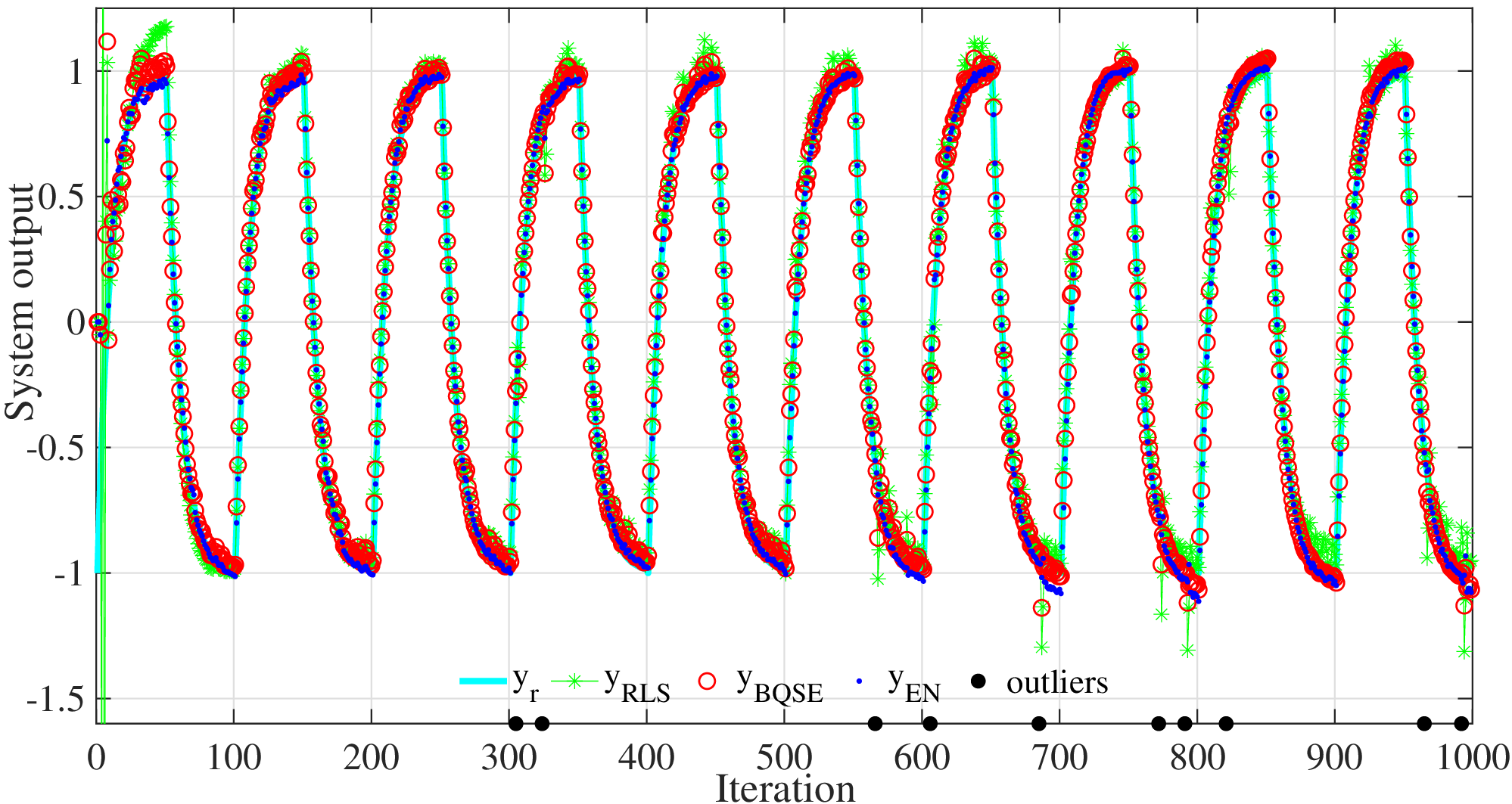}	
        \caption{The system output subject to Noise-\uppercase\expandafter{\romannumeral 3}}
	\label{fig8}	
\end{figure}

\begin{figure}[h]	
	\centering	
        \includegraphics[width=0.48\textwidth]{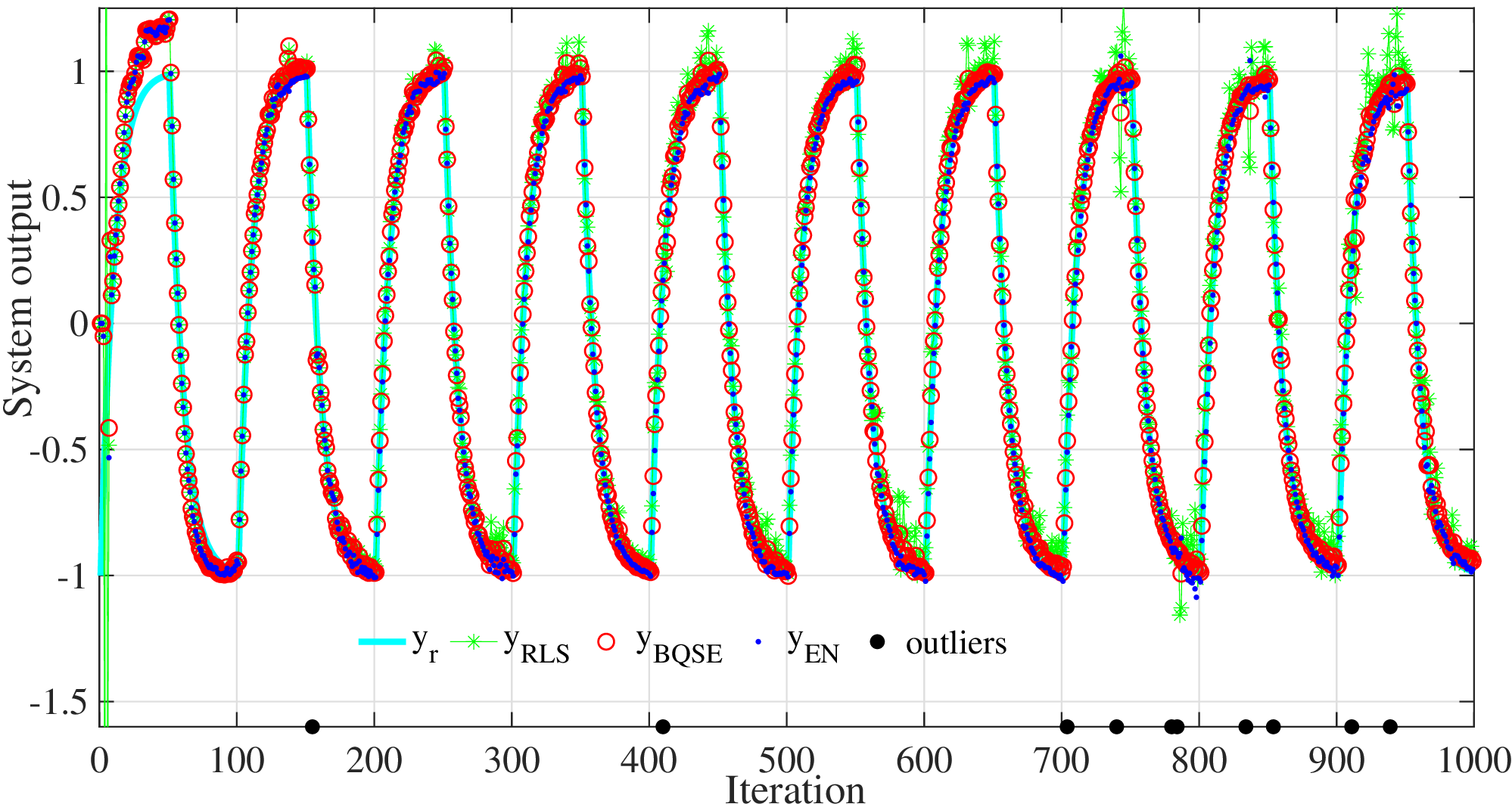}	
        \caption{The system output subject to Noise-\uppercase\expandafter{\romannumeral 4}}
	\label{fig9}	
\end{figure}

The system outputs by RLS, BQSE and the ensemble control subject to Noise-\uppercase\expandafter{\romannumeral 3} and Noise-\uppercase\expandafter{\romannumeral 4} are shown in Fig. \ref{fig8} and Fig. \ref{fig9}, respectively. We observe from Fig. \ref{fig8} that small spikes exist in the outputs by RLS and BQSE at the 566th, 685th, 772nd, 791st and 992nd iterations due to the existence of outliers. In Fig. \ref{fig9}, the outputs of RLS exhibit noticeable spikes at the 740th, 780th, 834th and 911st iterations with the presence of outliers, while the outputs of BQSE show relatively small spikes at the 155th, 740th and 834th iterations. In comparison, the proposed ensemble control renders no spike, demonstrating its robustness to outliers during the tracking control.




\section{Conclusions} \label{section5}

In this work, we proposed an adaptive ensemble control for stochastic systems with resilience against asymmetric noises and outliers. We constructed mixed ALD noise to model the distribution of noise, which can be generalized to represent symmetric and asymmetric noises and outliers. We designed an iterative quantile filter to learn the parameters of the system subject to mixed ALD noise, and derived an ensemble control law based on the iterative quantile filter. Numerical simulation results demonstrate that the proposed approach provides robustness against asymmetric noises and outliers for stochastic system control, under various tracking control scenarios and noise conditions. We plan to extend the proposed ensemble control with dual properties in our future work, aiming at reducing overshoots during rapid adaptation phases in the control. We also plan to integrate nonlinear system control in the ensemble controller in our extension.

\bibliography{mybibfileall}

\end{document}